\documentclass[12pt]{article}
\usepackage{amssymb,amsfonts,amsmath, psfrag,eepic}
\newtheorem{theo}{Theorem}

\newtheorem{coro}[theo]{Corollary}

\makeatletter \@addtoreset{equation}{section}
\@addtoreset{theo}{section} \makeatother

\def\qed{\hfill \rule{4pt}{7pt}}

\begin{document}
\begin{center}
{\large \bf Reduction of $m$-Regular Noncrossing Partitions}
\end{center}

\begin{center}
William Y. C. Chen$^1$, Eva Y. P. Deng$^2$ and Rosena R. X. Du$^3$
\end{center}
\vspace{0.2cm} \centerline{Center for Combinatorics, LPMC}
\centerline{ Nankai University, Tianjin  300071, P. R. China}
\vspace{3mm} \centerline{$^1$chen@nankai.edu.cn,
$^2$dengyp@eyou.com, $^3$du@nankai.edu.cn}

\vskip 3mm \centerline{Revised March 5, 2004}

\vspace{0.5cm}

\noindent {\bf Abstract.} In this paper, we present a reduction
algorithm which transforms  $m$-regular partitions of $[n]=\{1, 2,
\ldots, n\}$ to $(m-1)$-regular partitions of $[n-1]$. We show
that this algorithm preserves the noncrossing property. This
yields a simple explanation of an identity due to Simion-Ullman
and Klazar in connection with enumeration problems on noncrossing
partitions and RNA secondary structures. For ordinary noncrossing
partitions, the reduction algorithm leads to a representation of
noncrossing partitions in terms of independent arcs and loops, as
well as an identity of Simion and Ullman which expresses the
Narayana numbers in terms of the Catalan numbers.



\section{Introduction}

A {\it partition} $P$ of $[n]=\{1,2,...,n\}$ is a collection
$\{B_1, B_2, ..., B_k\}$ of nonempty disjoint subsets of $[n]$,
called {\it blocks} such that $B_1\cup \cdots \cup B_k=[n]$. We
may assume that $\{B_1, B_2, ..., B_k\}$ are listed in the
increasing order of their minimum elements. The set of all
partitions of $[n]$ with $k$ blocks is denoted by
$\mathcal{P}(n,k)$. The cardinality of $\mathcal{P}(n,k)$ is the
well-known Stirling number of the second kind \cite{ST1}.

A partition $P\in \mathcal{P}(n,k)$ is called $m$-{\it regular},
$m \geq 1$, if for any two distinct elements $x,y$ in the same
block, we have $|x-y| \geq m$. If all blocks of $P$ are singletons
(of cardinality one) we set $m=\infty$. The set of $m$-regular
partitions in $\mathcal{P}(n,k)$ is denoted by
$\mathcal{P}(n,k,m)$, and its cardinality  is denoted by
$p(n,k,m)$. When $m=1$,  a $1$-regular partition is an ordinary
partition.  A partition is called {\it poor} if each block
contains at most two elements. The set of all poor partitions in
$\mathcal{P}(n,k,m)$ is denoted by $\mathcal{P}_{2}(n,k,m)$, and
its cardinality is denoted by $p_{2}(n,k,m)$.

Any partition $P$ can be expressed by its {\it canonical
sequential form} $P=a_1a_2 \cdots a_n$ where $a_i=j$ if the
element $i$ is in the block $B_j$. For instance, $1231242$ is the
canonical sequential form of $P=(1,4)(2,5,7)(3)(6)\in
\mathcal{P}(7,4,2)$. In fact, one can use a sequence on any set of
$k$ symbols to represent a partition of $k$ blocks, where the
symbols are linearly ordered. If we use the alphabet $\{a,b,c,d\}$
of four letters with the order $a<b<c<d$, then the corresponding
canonical sequential form for $P$ becomes $abcabdb$. Note that if
$a_1a_2\cdots a_n$  is a canonical sequential form of a partition
with $k$ blocks, then each of $1, 2, \ldots, k$ appears at least
once and the first occurrence of $i$ precedes that of $j$ if
$i<j$. The sequence $a_1a_2\cdots a_n$ is also called the {\it
restricted growth function} of a partition $P$ \cite{WW}. These
two requirements  are the normalization conditions of the
Davenport-Schinzel sequences \cite{Dan65, MS}, as noted by Klazar
\cite{K1, K2}.

We say that $P \in \mathcal{P}(n,k,m)$ is $abab$-{\it free} if its
canonical sequential form does not contain any subsequence (not
necessarily a consecutive segment) of the form $\cdots a \cdots b
\cdots a \cdots b \cdots$, which is often written as $abab$.
Equivalently, $P$ is $abab$-free if there do not exist four
elements $x,y,u,v \in [n]$ with $x<u<y<v$ such that $x,y$ belong
to the same block and $u,v$ belong to another block. The set of
$abab$-free partitions in $ \mathcal{P}(n,k,m)$ is denoted by
$\mathcal{P}(abab;n,k,m)$. An $abab$-free partition is also called
a {\it noncrossing partition}. For example, the partition
$P=(1,4)(2,5,7)(3)(6)$ is not $abab$-free because of the violation
of the four elements $1<2<4<5$.

Regular partitions also arise in the enumeration of RNA secondary
structures. In biology, an RNA sequence can be viewed as a
sequence of molecules $A$ (adenine), $C$ (cytosine), $G$ (guanine)
and $U$ (uracil); these single-stranded molecules fold onto
themselves by the so-called {\it Watson-Crick rules}: $A$ forms
base pairs with $U$ and $C$ forms base pairs with $G$. A helical
structure can be formed based on the sequence of molecules and the
rules. If such a helical structure can be realized as a planar
graph, then it is called an {\it RNA secondary structure}. In the
mathematical modelling of RNA secondary structures, one may
disregard what the molecules are and consider a helical structure
as a sequence of numbers $1,2,...,n$ along with some base pairs,
where we have the restriction that all base pairs are allowed
except for any  two adjacent numbers $i$ and $i+1$, and there do
not exist two base pairs $(i,j)$ and $(k,l)$ with $i<k<j<l$. In
this setting, the set $\mathcal{R}(n,k)$ of all RNA secondary
structures with length $n$ and $k$ base pairs can be viewed as the
set $\mathcal{P}_{2}(abab;n,n-k,2)$ of noncrossing poor
partitions. A formula for $\mathcal{R}(n,k)$ is obtained by
Schmitt and Waterman \cite{SW1} in terms of the Narayana numbers.

However, a further biological consideration suggests a
generalization of the above mathematical modelling of RNA
secondary structures. As pointed out by Hofacker, Schuster and
Stadler \cite{HSS}, within each matching bracket (or base pair)
there should be at least three elements. In the language of
partitions, that is to say that if $i$ and $j$ are in the same
block, then we have $|i-j|\geq 4$. Equivalently, this is the
notion of $4$-regular partitions. Thus, we are led to the study of
$m$-regular noncrossing poor partitions. It is known that for
$m=1$, such partitions correspond to Motzkin paths, for $m=2$,
there is a correspondence with Motzkin paths without peaks
\cite{AN, VV}. In general, Klazar \cite{K1} gives a formula for
the number of $m$-regular noncrossing poor partitions.

The main result of this paper is a reduction algorithm that
transforms a partition in $\mathcal{P}(n,k,m)$ to a partition in
$\mathcal{P}(n-1,k-1,m-1)$. We show that the algorithm preserves
the noncrossing property or the $abab$-free property. This leads
to a quick explanation of the following identity:
\begin{equation}\label{mainr}
p(abab;n,k,m)=p_{2}(abab;n-1,k-1,m-1).
\end{equation}
An earlier version of this relation was first obtained by Simion
and Ullman \cite{SU1}, where the relation is stated for $m=2$.
Klazar found the above identity in general and gave a generating
function proof in \cite{K1}. Another bijective proof of
(\ref{mainr}) was found by Klazar \cite{K2}. We should note that
the notations in [4,5] are somewhat different. No simple
explanation of (\ref{mainr}) seems to be known. We hope that our
algorithm may have served this purpose.

It is worth noting that ordinary noncrossing partitions can be
further reduced into independent arcs and loops (defined
subsequently, just before Theorem \ref{m1}). Essentially, this
gives a correspondence between noncrossing partitions and
2-Motzkin paths, and leads to an identity expressing the Narayana
numbers in terms of the Catalan numbers due to Simion and Ullman
\cite{SU1}.

\section{The Reduction Algorithm}

We begin with a bijective understanding of an identity of Yang
\cite{Y} concerning the number of $m$-regular partitions of $[n]$.

\begin{theo} For $m\geq 2$, we have
\begin{equation}\label{m-m-1}
p(n,k,m) = p(n-1, k-1, m-1).
\end{equation}
\end{theo}

For the case $m=2$, a 2-regular partition is called a ``restricted
partition'' and (\ref{m-m-1}) was obtained by  Prodinger
\cite{PH}. Bijective proofs of (\ref{m-m-1}) for 2-regular
partitions are given by many people including Prodinger \cite{PH},
Yang \cite{Y}. However, these proofs do not seem to apply to
general $m$-regular partitions.

In this paper, we find a simple reduction algorithm for
$m$-regular partitions. The key idea is to use a digraph to
represent a partition, which is called {\it the linear
representation}. Given a partition $P=\{B_1, B_2, \ldots, B_k\}$
of $[n]$, we draw a digraph $D(P)$, or $D$ for short, on the
vertex set $[n]$. For each block $B_i$, we associate it with a
directed path $P_i$ starting with the minimum element in $B_i$,
and going through elements in $B_i$ in the increasing order. Note
that when a block $B_i$ has only one element, the corresponding
path is an isolated vertex. The digraph $D$ can be drawn on a line
such that the vertices $1, 2, \ldots, n$ are arranged in the
increasing order and the arcs always have the direction from left
to right. For this reason, one does not really need to display the
direction of each arc (see Figure 1). An undirected version of the
linear representation of a partition was used by Simion \cite{S1}.
As we will see, the directions are useful to clarify the argument
for the reduction algorithm.

\noindent{\bf The  Reduction Algorithm:} For a partition $P\in
\mathcal{P}(n,k,m)$, where $n,k,m\geq 1$, we may reduce it to a
partition in $\mathcal{P}(n-1, k-1, m-1)$:
\begin{itemize}
\item[1.] For each arc $(i,j)$ in the linear representation of
$P$, replace it by the arc $(i,j-1)$.
\item[2.] Delete the vertex $n$.
\end{itemize}

\begin{theo}\label{biject}
When $m \geq 2$, the reduction algorithm gives a bijection between
$\mathcal{P}(n,k,m)$ and $\mathcal{P}(n-1,k-1,m-1)$.
\end{theo}

\noindent{\it Proof.} Suppose that $P$ is a partition in
$\mathcal{P}(n,k,m)$. Let $D$ be the linear representation of $P$,
and let $D^{\prime}$ be the digraph obtained from $D$ by reducing
every arc (replacing $(i,j)$ by $(i, j-1)$). Since $m\geq 2$, it
is clear that in $D^{\prime}$ every arc has the direction from the
smaller vertex to the bigger vertex, and for each vertex $j$ in
$D^{\prime}$, neither its indegree nor outdegree is greater than
1. Thus, each component of $D^{\prime}$ is a directed path from
the minimum vertex to the maximum vertex in the increasing order.
In other words, $D^{\prime}$ is also a linear representation of an
$(m-1)$-regular partition of $[n]$. It is easy to see that
$D^{\prime}$ has the same number of arcs as $D$. It follows that
$D^{\prime}$ and $D$ have the same number of connected components.
Let $H$ be the digraph obtained from $D^{\prime}$ by deleting the
isolated vertex $n$.  Then $H$ is the linear representation of the
desired partition.

By reversing the above procedure, one may show that the reduction
algorithm yields a bijection. \qed

An example is given in Figure 1, which illustrates the bijection
between $\mathcal{P}(5,3,2)$ and $\mathcal{P}(4,2,1)$.
\begin{figure}[h,t]
\begin{center}
\begin{picture}(360,220)
\setlength{\unitlength}{7mm} \thicklines
\put(-0.2,-0.3){\line(1,0){17}} \put(-0.2,10){\line(1,0){17}}
\put(-0.2,11){\line(1,0){17}} \put(-0.2,-0.3){\line(0,1){11.3}}
\put(9,-0.3){\line(0,1){11.3}} \put(16.8,-0.3){\line(0,1){11.3}}
\put(3.5,10.3){$\mathcal{P}(5,3,2)$}
\put(12,10.3){$\mathcal{P}(4,2,1)$} \thinlines
\put(3.5,-0.3){\line(0,1){10.3}} \put(12.5,-0.3){\line(0,1){10.3}}

\multiput(4,0)(1,0){5}{\circle*{0.1}} \qbezier(4,0)(6,2)(8,0)
\qbezier(5,0)(6,1)(7,0) \put(0,0){(1,5)(2,4)(3)}
\multiput(13,0)(1,0){4}{\circle*{0.1}}
\qbezier(13,0)(14.5,1.5)(16,0) \qbezier(14,0)(14.5,0.5)(15,0)
\put(9.5,0){(1,4)(2,3)}

\multiput(4,1.5)(1,0){5}{\circle*{0.1}}
\qbezier(5,1.5)(6,2.5)(7,1.5) \qbezier(6,1.5)(7,2.5)(8,1.5)
\put(0,1.5){(1)(2,4)(3,5)}
\multiput(13,1.5)(1,0){4}{\circle*{0.1}}
\qbezier(14,1.5)(14.5,2.1)(15,1.5)
\qbezier(15,1.5)(15.5,2.1)(16,1.5) \put(9.5,1.5){(1)(2,3,4)}

\multiput(4,3)(1,0){5}{\circle*{0.1}} \qbezier(4,3)(5.5,4.3)(7,3)
\qbezier(6,3)(7,3.8)(8,3) \put(0,3){(1,4)(2)(3,5)}
\multiput(13,3)(1,0){4}{\circle*{0.1}} \qbezier(13,3)(14,4)(15,3)
\qbezier(15,3)(15.5,3.5)(16,3) \put(9.5,3){(1,3,4)(2)}

\multiput(4,4.5)(1,0){5}{\circle*{0.1}}
\qbezier(4,4.5)(5.5,5.7)(7,4.5) \qbezier(5,4.5)(6.5,5.7)(8,4.5)
\put(0,4.5){(1,4)(2,5)(3)}
\multiput(13,4.5)(1,0){4}{\circle*{0.1}}
\qbezier(13,4.5)(14,5.5)(15,4.5) \qbezier(14,4.5)(15,5.5)(16,4.5)
\put(9.5,4.5){(1,3)(2,4)}

\multiput(4,6)(1,0){5}{\circle*{0.1}} \qbezier(4,6)(5,7)(6,6)
\qbezier(6,6)(7,7)(8,6) \put(0,6){(1,3,5)(2)(4)}
\multiput(13,6)(1,0){4}{\circle*{0.1}}
\qbezier(13,6)(13.5,6.6)(14,6) \qbezier(15,6)(15.5,6.6)(16,6)
\put(9.5,6){(1,2)(3,4)}

\multiput(4,7.5)(1,0){5}{\circle*{0.1}}
\qbezier(4,7.5)(5,8.5)(6,7.5) \qbezier(5,7.5)(6.5,8.5)(8,7.5)
\put(0,7.5){(1,3)(2,5)(4)}
\multiput(13,7.5)(1,0){4}{\circle*{0.1}}
\qbezier(13,7.5)(13.5,8.1)(14,7.5)
\qbezier(14,7.5)(15,8.5)(16,7.5) \put(9.5,7.5){(1,2,4)(3)}

\multiput(4,9)(1,0){5}{\circle*{0.1}} \qbezier(4,9)(5,10)(6,9)
\qbezier(5,9)(6,10)(7,9) \put(0,9){(1,3)(2,4)(5)}
\multiput(13,9)(1,0){4}{\circle*{0.1}}
\qbezier(13,9)(13.5,9.6)(14,9) \qbezier(14,9)(14.5,9.6)(15,9)
\put(9.5,9){(1,2,3)(4)}
\end{picture}
\end{center}
\caption{Correspondence between $\mathcal{P}(5,3,2)$ and
$\mathcal{P}(4,2,1)$. }
\end{figure}

\section{Reduction of Noncrossing Partitions }

In this section, we show that the reduction algorithm preserves
the noncrossing property. This gives a simple explanation of the
following identity due to Simion and Ullman \cite{SU1} and Klazar
\cite{K1}.

\begin{theo} \label{NC_poor}
For $m\geq 2$, we have
\begin{equation}\label{abab}
p(abab;n,k,m)=p_{2}(abab;n-1,k-1,m-1).
\end{equation}
\end{theo}

\noindent{\it Proof.} Suppose $P$ is a partition in
$\mathcal{P}(abab;n,k,m)$ and $D$ is the linear representation of
$P$. Let $Q$ be the partition of $[n-1]$ obtained from $P$ by
applying the reduction algorithm, and  $D^{\prime}$  the linear
representation of $Q$. Suppose $D^{\prime}$ has a path of length
two, $i\rightarrow j \rightarrow k$, say, then $(i, j+1)$ and $(j,
k+1)$ are two arcs in $D$. Since $D$ is a linear representation,
these two arcs $(i, j+1)$ and $(j, k+1)$ belong to different
components, which contradicts the assumption that $P$ is
$abab$-free. It follows that $Q$ is a poor partition.

By the reduction algorithm we see that $Q$ has $k-1$ blocks and is
$(m-1)$-regular. It remains to show that $Q$ is noncrossing.
Suppose that there are four elements $x< u< y < v$ such that
$B_i=\{x,y\}$ and $B_j=\{u,v\}$, where $B_i$ and $B_j$ are
different blocks of $Q$. Then in the linear representation
$D^{\prime}$, $(x, y)$ and $(u, v)$ are two crossing arcs, it
follows that $(x, y+1)$ and $(u, v+1)$ are two crossing arcs in
$D$, which is a contradiction to the assumption that $P$ is
noncrossing. The converse can be justified in the same manner.
Therefore, we have established the desired one-to-one
correspondence. \qed

In Figure 1, there are only two partitions in
$\mathcal{P}(abab;5,3,2)$: $(1,3,5)(2)(4)$ and $(1,5)(2,4)(3)$;
after applying the reduction algorithm, we get $(1,2)(3,4)$ and
$(1,4)(2,3)$, which are the two partitions in
$\mathcal{P}_{2}(abab;4,2,1)$.

Theorem \ref{NC_poor} is useful for the enumeration of $m$-regular
noncrossing partitions \cite{K1}. We recall that the notion of
$m$-regular noncrossing poor partitions coincides with that of
general RNA secondary structures. The reduction algorithm can be
used even for ordinary noncrossing partitions. In a digraph $D$,
we say that two arcs are {\it independent} if they have no vertex
in common, and a {\it loop} is an arc from a vertex to itself.

\begin{theo} \label{m1}
There is a one-to-one correspondence between noncrossing
partitions of $[n]$ with $k$ blocks  and  digraphs on $[n-1]$
consisting of $n-k$ independent noncrossing arcs or loops.
\end{theo}

Digraphs described in Theorem \ref{m1} are related to $2$-Motzkin
paths introduced by Barcucci, del Lungo, Pergola and Pinzani
\cite{EB1}. Roughly speaking, if we consider loops and singletons
as  straight and wavy level steps respectively, then we obtain
\begin{theo}\label{theo33}
There is a one-to-one correspondence between noncrossing
partitions of $[n]$ with $k$ blocks  and  $2$-Motzkin paths of
length $n-1$ with $n-k$ straight level steps or up steps.
\end{theo}

Figure 2 is an illustration of the above bijection.
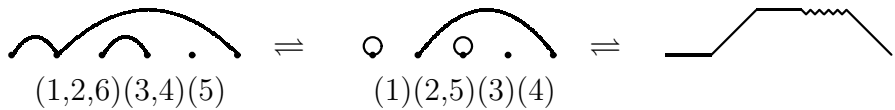
\begin{figure}[h,t]
\begin{center}
\begin{picture}(320,20)
\setlength{\unitlength}{6mm} \thicklines
\multiput(0,0)(1,0){6}{\circle*{0.1}} \qbezier(0,0)(0.5,0.8)(1,0)
\qbezier(1,0)(3,2)(5,0) \qbezier(2,0)(2.5,0.8)(3,0)
\put(0.5,-1){(1,2,6)(3,4)(5)} \put(5.8,0){$\rightleftharpoons$}
\multiput(8,0)(1,0){5}{\circle*{0.1}} \put(8,0.2){\circle{0.4}}
\put(10,0.2){\circle{0.4}} \qbezier(9,0)(10.5,2)(12,0)
\put(8.0,-1){(1)(2,5)(3)(4)} \put(12.8,0){$\rightleftharpoons$}
\put(14.5,0){\line(1,0){1}} \put(15.5,0){\line(1,1){1}}
\put(16.5,1){\line(1,0){1}}
\multiput(17.6,0.9)(0.2,0){5}{\line(1,1){0.1}}
\multiput(17.5,1)(0.2,0){5}{\line(1,-1){0.1}}
\put(18.5,1){\line(1,-1){1}}
\end{picture}
\end{center}
\caption{A noncrossing partition and the corresponding $2$-Motzkin
path.}
\end{figure}

In \cite{DS1}, Deutsch and Shapiro established a bijection between
ordered trees and $2$-Motzkin paths, and derived many important
consequences regarding combinatorial structures such as Dyck
paths, bushes, $\{0,1,2\}$-trees, Schr\"oder paths, RNA secondary
structures, noncrossing partitions, Fine paths, etc. The above
theorem on the reduction of noncrossing partitions to 2-Motzkin
paths can be viewed as a simpler version of the Deutsch-Shapiro
correspondence, since there are easy bijections between ordered
trees and noncrossing partitions \cite{ZAKs, PH2}.

Theorem \ref{theo33} leads to an identity of Simion and Ullman
\cite{SU1} expressing the Narayana numbers by the Catalan numbers.
Recall that the Catalan number $C_n = {1 \over n+1} \, {2n \choose
n}$ counts the number of plane trees with $n+1$ vertices, and the
Narayana number $N_{n,k} = {1 \over n} \, {n\choose k} \, {
n\choose k-1}$ is the number of plane trees with $n+1$ vertices
and $k$ leaves, which also counts the number of noncrossing
partitions of $[n]$ with $k$ blocks \cite{ZAKs, PH2}.

\begin{coro}[Simion and Ullman \cite{SU1}, Corollary 3.2] For all
$n\geq 1$ and $1 \leq k \leq n$, we have the following relation:
\begin{equation}\label{n-c}
N_{n,k} = \sum_{i=0}^{n-k} {n-1 \choose 2i} {n-2i-1 \choose n-i-k}
C_{i}.
\end{equation}
\end{coro}

\noindent{\it Proof.} Suppose $P$ is a noncrossing partition on
$[n]$ with $k$ blocks. Let $H$ be the digraph on $[n-1]$ with
independent noncrossing arcs, namely the linear representation of
the partition obtained from $P$ by applying the reduction
algorithm. Suppose $i$ is the number of loops in $H$. Removing the
loops, we get a digraph $H^{\prime}$ on $n-i-1$ vertices with
$k-1$ components consisting of $n-k-i$ independent arcs and
$2k+i-n-1$ isolated vertices.  The digraph $H^{\prime}$
corresponds to a noncrossing poor partition on $[n-i-1]$ with
$k-1$ blocks (namely $2k+i-n-1$ singletons and $n-k-i$ blocks of
size two). Hence we get
\begin{equation}
N_{n,k} = \sum_{i=0}^{n-k}\, {n-1 \choose i} \,
                          p_{2} (abab; n-i-1,k-1,1).
\end{equation}
There is a bijection between noncrossing poor partitions of
$[2n-2k-2i]$ without singletons and Dyck paths of length
$2n-2k-2i$ (see \cite{ST2}, p. 222, Exercise 6.19 (n) and its
solution on p. 258), and it is well-known that the number of such
Dyck paths equals the $(n-k-i)$-th Catalan number. In view of the
number of ways to choose the singletons, we obtain
\[ p_{2} (abab; n-i-1,k-1,1)= {n-i-1 \choose 2n-2i-2k} C_{n-i-k}.\]
By replacing $n-i-k$ with $i$ and taking summation,  we get
(\ref{n-c}). \qed

\vspace{0.5cm} \noindent{\bf Acknowledgments.} The authors would
like to thank E. Deutsch, M. Klazar, and J. Zeng for helpful
comments. We also thank the referees for important suggestions
adopted in the revised version. This work was done under the
auspices of the 973 Project on Mathematical Mechanization, and the
National Science Foundation of China.

\end{document}